\pdfoutput=1

\documentclass[letterpaper]{article}

\usepackage{arxiv}
\usepackage{amssymb}
\usepackage{amsmath}
\usepackage{bm}
\usepackage{amsthm}
\usepackage{mathtools}
\usepackage{natbib}
\usepackage{extarrows}
\usepackage{float}
\usepackage{array}
\usepackage{hyperref}
\usepackage{booktabs}
\usepackage{multirow}
\usepackage{siunitx}
\usepackage{caption}

\title{Autocovariance Estimation in the Presence of Changepoints}

\author{
  Colin Gallagher\\
  Department of Statistics and Operation Research\\
  Clemson University\\
   \And
 Rebecca Killick \\
  Department of Statistics\\
  Lancaster University\\
     \And
 Robert Lund \\
  Department of Statistics\\
  University of California, Santa Cruz\\
       \And
 Xueheng Shi \\
  Department of Statistics\\
  University of California, Santa Cruz\\
  \texttt{xshi38@ucsc.edu}  
}

\newtheorem{thm}{Theorem}

\newtheorem{col}{Corollary}

\DeclareMathOperator*{\median}{median}
\newcommand\ddfrac[2]{\frac{\displaystyle #1}{\displaystyle #2}}

\begin{document}
\maketitle

\begin{abstract}
This article studies estimation of a stationary autocovariance structure in the presence of an unknown number of mean shifts.  Here, a Yule-Walker moment estimator for the autoregressive parameters in a dependent time series contaminated by mean shift changepoints is proposed and studied. The estimator is based on first order differences of the series and is proven consistent and asymptotically normal when the number of changepoints $m$ and the series length $N$ satisfies $m/N \rightarrow 0$ as $N \rightarrow \infty$.  
\end{abstract}


\section{Introduction}

Time series dynamics often change due to external events or internal systematic fluctuations.  Changepoint analyses have become commonplace in identifying whether and when abrupt changes take place.

Evolving from the original treatise for a single location parameter shift in \cite{Page1954}, the majority (but not all) changepoint analyses check for shifts in the series' mean.  Since then, considerable work on changepoints has been conducted, including applications to climatatology \citep{lund2017temp}, economics \citep{NorwoodKillick2018}, and disease modelling \citep{hall2000disease}. This paper studies approaches for estimating an unknown autocorrelation function when an unknown number of mean shift changepoints are present.  

Many changepoint techniques assume independent and identically distributed (IID) model errors; however, time series data are typically correlated (e.g., daily temperaturs and stock prices and DNA sequences). Changepoint techniques may incorrectly segment the data (estimate too many or too few changepoints) should non-zero autocorrelation be ignored. Two primary ways of tackling this issue have been pursued: 1) including autocorrelation in the multiple changepoint model \citep{lilund2012MDL}, and 2) decorrelating the time series prior to any changepoint analysis \citep{chakar2017, shi2021comparison}.  In either case, one needs the autocovariance function of the data.  With a good estimate of the series' autocovariance, one-step-ahead prediction residuals can be estimated --- and these residuals are always uncorrelated (independent up to estimation for Gaussian series). Indeed, a principle of \citep{shi2021comparison, robbins2011} is that good multiple changepoint detection routines can be devised by applying IID methods to one-step-ahead prediction residuals.

Some multiple changepoint models for time series allow all model parameters, including those governing the correlation structure of the series, to change at each changepoint time.  These scenarios are easier to handle (computationally quicker) as dynamic programming techniques can quickly optimize penalized objective functions; see \citet{killick2012pelt} and \citet{maidstone2017fpop}. The key here is that the objective function that is optimized needs to be additive in its segments (regimes). 

A more parsimonious setup allows means to shift with each changepoint time, but keeps error autocovariances constant across all regimes.  As this breaks the objective function additivity required, fast dynamic programming techniques cannot be directly applied.  Thus, estimates of global autocovariances governing all regimes need to be developed.  Standard estimation schemes degrade unless the changepoint locations are a priori known (known changepoint times allow the segment means to be better estimated and subtracted from the series a priori to autocovariance calculation).  Should process autocovariance estimators be off (due to bias or inconsistency), then decorrelation techniques will also suffer.  As such, there is a need to estimate a stationary autocovariance function in the presence of mean shifts.

This paper studies autocovariance estimation in the presence of mean shifts in more detail.  We devise a method based on first order differencing that beats rolling window methods.  The scenario is asymptotically quantified when model errors are drawn from a causal autoregressive process.   The rest of this paper proceeds as follows.  The next section narrates our setup and discusses approaches to the problem.  Section 3 then develops an estimation technique based on lag one differences of the series. Section 4 proves consistency and asymptotic normality of these estimators and Section 5 assesses their performance in simulations.  Section 6 applies the results to an annual precipitation series and Section 7 concludes the paper with comments.

\section{Model and Estimation Approaches}

Suppose that $\{ X_t \}_{t=1}^N$ is a time series having an unknown number of mean shift changepoints, denoted by $m$, occurring at the unknown ordered times $1 < \tau_1 < \tau_2 < \cdots < \tau_m \le N$. These $m$ changepoints partition the series into $m+1$ distinct segments, each segment having its own distinct mean. The model can be written as 
\begin{equation}
\label{shiftmodel}
X_t=\kappa_{s(t)} + \epsilon_t.
\end{equation}
If $s(t)$ denotes the series' segment at time $t$, which takes on values in $\{0, 1, \cdots, m \}$, then $\kappa_{s(t)} = \mu_i$ is constant for all times in the $i^{th}$ regime:
\begin{align*}
\kappa_{s(t)}=
     \begin{cases}
     \;\mu_0, \quad &1 \le t \le \tau_1,\\
     \;\mu_1, \quad &\tau_1+1 \le t \le \tau_2,\\
     \; \; \qquad \vdots  \\
     \; \mu_m, \quad &\tau_{m}+1 \le t \le N.
     \end{cases}
\end{align*} 
We assume that $\{ \epsilon_t \}$ is a stationary causal AR$(p)$ time series that applies to all regimes.  Then $\{ \epsilon_t \}$ obeys  
\begin{align}
    \epsilon_t=\phi_1 \epsilon_{t-1}+\cdots +\phi_{p}\epsilon_{t-p}+Z_t,
    \quad t \in \mathbb{Z},
\label{ARp}
\end{align}
where $\{ Z_t \}$ is IID white noise with zero mean, variance $\sigma^2$, and a finite fourth moment (this enables consistent estimation of the autoregressive parameters $\phi_1, \ldots , \phi_p$).  While more general ARMA($p,q$) $\{ \epsilon_t \}$ could be considered, we work with AR($p$) errors because this model class is dense in all stationary short-memory series and estimation, prediction, and forecasting are easily conducted.  Adding a moving-average component $q \geq 1$ induces considerably more work and is less commonly found in changepoint applications.  
Few authors have considered autocovariance estimation in this changepoint corrupted setting. \cite{chakar2017}, in studying the AR(1) case, develops an estimator that is robust to mean shifts:
\begin{equation}
  \hat{\phi} = \ddfrac{\displaystyle \left ( \median_{1 \le t \le N-2} |X_{t+2} -X_t| \right)^2} {\displaystyle \left(\median_{1 \le t \le N-1} |X_{t+1} -X_t| \right)^2} - 1.
  \label{ar1robust}
\end{equation}
This is claimed to be consistent and satisfy the central limit theorem. It is not clear how to extend this work to cases where $p > 1$; moreover, we show that this estimator does not work particularly well in AR$(1)$ settings.

A second way to handle the problem jointly estimates the AR model and the number and locations of the changepoints \citep{lilund2012MDL}. For one way to do this, a penalized likelihood function can be optimized over all admissible changepoint configurations.  If the changepoint configuration is specified (known), the sequence is easily centered to a zero mean by subtracting segment sample means.  Any ARMA$(p,q)$ error model can be fitted to the centered series.  A genetic algorithm (GA) can be used to estimate the optimal penalized likelihood over all changepoint configurations under a variety of penalization schemes \citep{lund2020wbs}.  Unfortunately, the GA may take significant computing time to find the optimum:  there are $2^{N-1}$ different multiple changepoint configurations in a series of length $N$ to search over.  Sometimes, the GA may fail to converge at all. 

A third approach employs rolling windows. See \cite{marriott1954} and \cite{orcutt1969} for estimation of autocorrelations via window methods and \cite{BeaulieuKillick2018} for applications of these methods to changepoint problems in climatology.  For the window length $w$, with $w \le N$, a moving window scheme generates $N-w+1$ subsegments, the $i^{th}$ subsegment containing the points at the times $i, \ldots, i+w-1$.  Each subsegment is treated as a stationary time series (even though some may contain mean shifts and are thus truly nonstationary) and time series parameters are estimated in subsegment $i$ via the data in this subsegment only. The final estimates are taken as medians of the estimates over all subsegments. The hope is that most windows will be changepoint free, and medians over all subsegments will not be heavily influenced by the few windows that contain changepoints.  Of course, such a scheme may not use all data efficiently in estimation.  Moreover, \cite{BeaulieuKillick2018} demonstrated that the success of this estimation depends heavily on the choice of $w$.

\section{Moment Estimates based on Differencing}
\begin{center}
\begin{figure}
\centering
\includegraphics[scale=0.5]{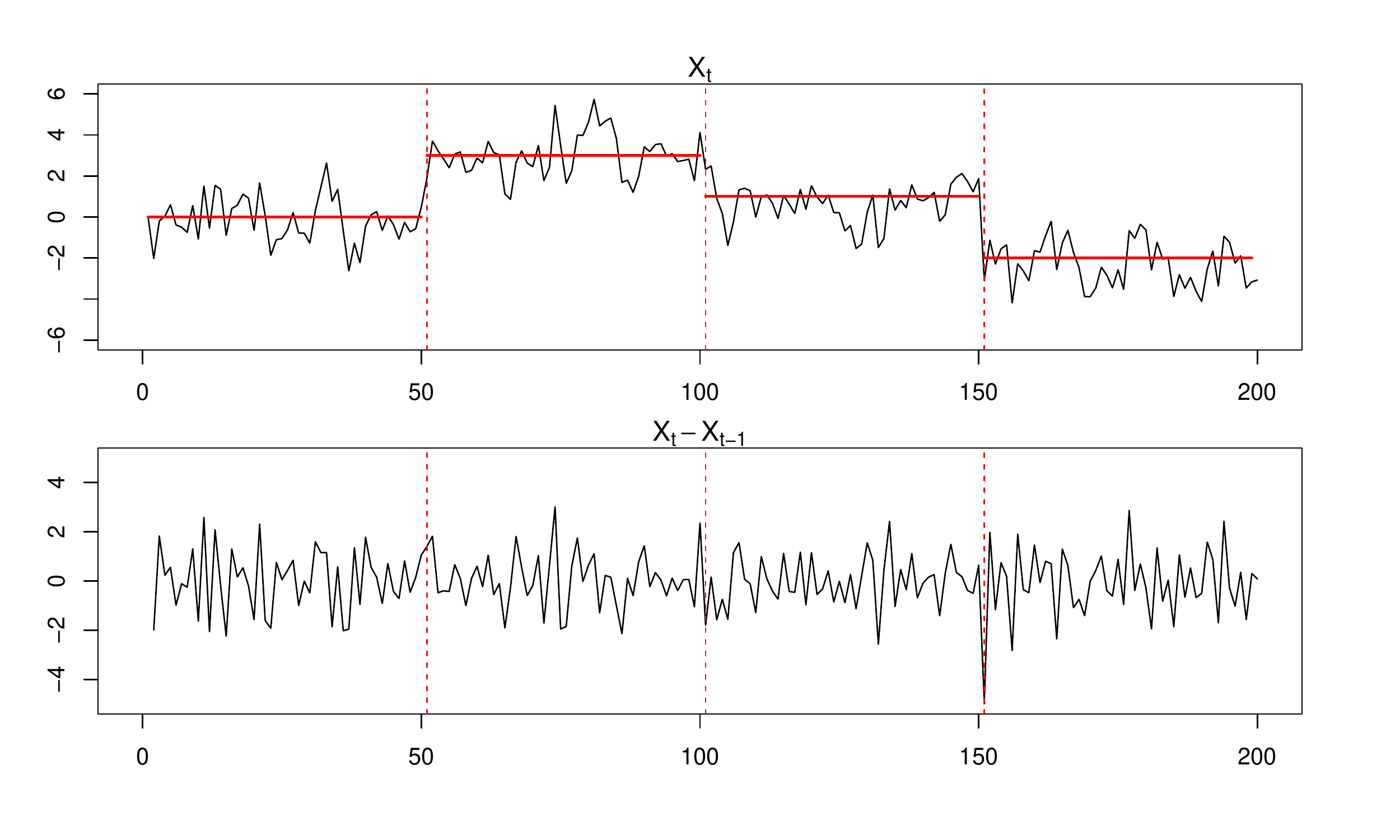}
\caption{An AR$(1)$ series $\{ X_t \}$ with three changepoints at the times $t=51, 101, 151$ (top panel).  Its first order difference (bottom panel) contains only three points with a non zero mean.}
\label{diff_demo}
\end{figure}
\end{center}

This section derives a system of linear equations that relates the autocorrelations of the differenced series to the AR($p$) coefficients.  First order differencing a series eliminates any piece-wise constant mean except at times where shifts occur; Figure \ref{diff_demo} illustrates the simple idea. Authors have previously used differencing to estimate global parameters in the changepoint literature.  For examples, \citet{fryzlewicz2014wbs} uses differencing to get an estimate of $\mbox{Var}(X_t)$, although IID errors are assumed in this work; (\ref{ar1robust}) from \cite{chakar2017} is also based on differencing.  However, there seems to be no previous literature using differencing to estimate AR($p$) parameters in a setting corrupted by mean shifts. As an aside, differencing also detrends a time series; the estimators below perform well if a time series contains both changepoints and a linear trend. 

A system of Yule-Walker moment equations can be used to estimate the AR($p$) parameters.  Our equations have the form 
\begin{equation}
\label{YULE}
\boldsymbol{\rho}_d=\boldsymbol{M}\boldsymbol{\phi},
\end{equation}
where $\boldsymbol{\rho}_d=[\rho_d(1)+1/2, \rho_d(2), \rho_d(3), \ldots, \rho_d(p)]^T$, $\rho_d(h)$ is the autocorrelation of the differenced series at lag $h$, $\boldsymbol{\phi}= [\phi_1, \phi_2, \phi_3, \cdots, \phi_p]^T$, and 
\begin{equation}
\label{Meq}    
\mathbf{M}=
\begin{bmatrix}
\frac{1}{2} &-\frac{1}{2} &-\left(\frac{1}{2}+\rho_d(1)\right) &\cdots &-\left(\frac{1}{2}+\sum_{j=1}^{p-2}\rho_d(j)\right) \\
\rho_d(1) &\rho_d(0) &\rho_d(1) &\cdots &\rho_d(p-2)\\
\rho_d(2) &\rho_d(1) &\rho_d(0) &\cdots &\rho_d(p-3)\\
\vdots &\vdots  & \vdots &\ddots &\vdots\\
\rho_d(p-1) &\rho_d(p-2) &\rho_d(p-3) &\cdots &\rho_d(0)\\
\end{bmatrix} 
\end{equation}  
This system is shown below to have a unique solution. 

To derive the system of equations in \ref{YULE}, make the notations $\gamma_X(h)=\mbox{Cov}(X_t, X_{t-h})$, $d_t = \epsilon_t - \epsilon_{t-1}$, and $\gamma_d(h) = \mbox{Cov}(d_t, d_{t-h})$. The covariance of $\{ d_t \}$ and $\{ X_t \}$ are related via
\[
\gamma_d(h) = 2 \gamma_X(h) - \gamma_X(h+1)-\gamma_X(h-1).
\]
Define $r(k)=\gamma_X(k)-\gamma_X(k-1)$.  Then $r(1)= -\gamma_d(0)/2$ and $r(k)=r(k-1)-\gamma_d(k-1)$ for $k \geq 2$.  Solving this difference equation recursively yields
\[
r(k) = -\frac{\gamma_d(0)}{2}-\sum_{\ell=1}^{k-1} \gamma_d(\ell),
\]
which implies that
\begin{equation}
\label{step1}
    \frac{r(k)}{\gamma_d(0)}=
    - \left[ \frac{1}{2}+ \sum_{\ell=1}^{k-1} \rho_d(\ell) \right].
\end{equation}
Combining with the AR$(p)$ definition in (\ref{ARp}), we extract
\begin{align*}
(a) \qquad &\gamma_X(h)=\phi_1\gamma_X(h-1)+\phi_2\gamma_X(h-2)+\cdots+\phi_p\gamma_X(h-p),\qquad &h=1,\ldots, p;\\
(b) \qquad &d_t=\phi_1 d_{t-1}+\cdots+\phi_p d_{t-p}+Z_t-Z_{t-1}, \qquad &t=0,\pm 1,\ldots;\\
(c) \qquad &\gamma_d(h)=\phi_1\gamma_d(h-1)+\phi_2\gamma_d(h-2)+\cdots+\phi_p\gamma_d(h-p),\qquad &h=2,\ldots, p;\\
\end{align*} 
where (c) follows by multiplying both sides of (b) by $d_{t-h}$ and taking expectations; (b) establishes $\{ d_t \}$ as a non-invertible ARMA($p,1)$ sequence.  

Subtracting (a) with $h=2$ from (a) with $h=1$ produces
\[
\gamma_X(1)-\gamma_X(2)=\phi_1[\gamma_X(0)-\gamma_X(1)]+\phi_2[\gamma_X(1)-\gamma_X(0)]+ \cdots+ \phi_p[\gamma_X(p-1)-\gamma_X(p-2)].
\]
Dividing this relation by $\gamma_d(0)$ and using (\ref{step1}) gives
\begin{equation}
\label{M1}
    \rho_d(1)+\frac{1}{2}=\frac{\phi_1}{2}-\frac{\phi_2}{2}-\phi_3\left[\frac{1}{2}+\rho_d(1)\right]-\cdots-\phi_p \left[\frac{1}{2}+\rho_d(1)+\cdots+\rho_d(p-2)\right],
\end{equation}
where $\rho_d(h)=\gamma_d(h) / \gamma_d(0)$ denotes the lag $h$ autocorrelation of $\{ d_t \}$.  To get another $p-1$ moment equations, divide (c) by $\gamma_d(0)$:
\begin{equation}
\label{M2}
\left(\begin{array}{c} \rho_d(2)\\ \vdots \\ \rho_d(p) \end{array}\right)=
\phi_1 \left(\begin{array}{c} \rho_d(1)\\ \vdots\\ \rho_d(p-1)\end{array}\right)
+A \left(\begin{array}{c} \phi_2 \\ \vdots \\ \phi_p \end{array} \right),
\end{equation}
where
\[
\mathbf{A}=
\begin{bmatrix}
\rho_d(0) &\rho_d(1) &\cdots &\rho_d(p-2)\\
\rho_d(1) &\rho_d(0) &\cdots &\rho_d(p-3)\\
\vdots  & \vdots &\ddots &\vdots\\
\rho_d(p-2) &\rho_d(p-3) &\cdots &\rho_d(0)\\
\end{bmatrix}. 
\]
Equation (\ref{YULE}) simply writes (\ref{M1}) and (\ref{M2}) into a single linear system.  

We now show that (\ref{M1}) and (\ref{M2}) provide a unique solution for the AR coefficients in terms of the autocorrelations of the differenced series. To do this, note that $\mathbf{A}$ is invertible, implying that
\begin{equation}
\label{sol1}
\left(\begin{array}{c} \phi_2 \\ \vdots \\ \phi_p \end{array} \right)=
\mathbf{A}^{-1}\left[ \left(\begin{array}{c} \rho_d(2)\\ \vdots \\ \rho_d(p) \end{array}\right)-\phi_1\left(\begin{array}{c} \rho_d(1)\\ \vdots\\ \rho_d(p-1)\end{array}\right) \right].
\end{equation}
To combine this with (\ref{M1}), define the vector
\[
\mathbf{c}^T=\left(-\frac{1}{2}, -\left(\frac{1}{2}+\rho_d(1)\right), \cdots, -\left(\frac{1}{2}+\sum_{j=1}^{p-2}\rho_d(j)\right)\right),
\]
and observe that
\begin{equation}
\label{sol2} 
\phi_1 \left[  
\frac{1}{2}+\mathbf{c}^T\mathbf{A}^{-1}
\left(
\begin{array}{c} 
\rho_d(1) \\ 
\vdots    \\ 
\rho_d(p-1)
\end{array}
\right) \right] =
\rho_d(1) +\frac{1}{2} + \mathbf{c}^T \mathbf{A}^{-1}
\left(
\begin{array}{c} 
\rho_d(2)\\ 
\vdots \\ 
\rho_d(p)
\end{array}
\right).
\end{equation}
We see that given the first $p$ autocorrelations of the differenced series, there is a unique set of parameters $\phi_1, \ldots, \phi_p$ satisfying  (\ref{sol1}) and (\ref{sol2}), and thus satisfying (\ref{YULE}).

To estimate the AR($p$) parameters in practice, the sample autocorrelations of $\{ d_t \}$ are calculated for the lags $0, 1, \ldots , p$.  Our estimator has the form 
\begin{equation}
\label{YWestimator}
\hat{\mathbf{\phi}}=\hat{\mathbf{M}}^{-1}\hat{\mathbf{\rho}}_d,
\end{equation}
where the elements of $\hat{\boldsymbol{M}}$ and $\hat{\boldsymbol{\rho}}_d$ simply replace $\rho_d(h)$ with its estimator
\[
\hat{\rho}_d(h) = 
\frac{\hat{\gamma}_d(h)}{\hat{\gamma}_d(0)}=
\frac
{\sum_{t=2}^{N-h} (d_t-\bar{d})(d_{t+h}-\bar{d})}
{\sum_{t=1}^N (d_t-\bar{d})^2}.
\]
Here, $\bar{d}= (N-1)^{-1}\sum_{t=1}^{N-1} d_t$ is the sample mean, which is included here to allow for the possibility of a linear trend in the original series. 

Invertibility of $\hat{\mathbf{M}}$ follows an argument nearly identical to that for $\mathbf{M}$, and is based on the fact that the matrix of sample autocorrelations $\hat{\mathbf{A}}$ is invertible with probability one (see Chapter 8 of \cite{BrockwellDavis1991}). 

If the number of changepoints $m$ is small relative to $N$, then the mean shifts will have a negligible impact on the estimated covariance of the differences, since $X_t-X_{t-1} = d_t-d_{t-1}$ except at the changepoint times $\tau_1, \ldots, \tau_m$.

We end this section by estimating $\sigma^2$.  Multiplying both sides of (b) by $d_{t-1}$, taking expectations, and solving for $\mbox{Var}(Z_t) = \sigma^2$, yields
\[
\sigma^2=\left(\sum_{j=1}^p \phi_j \gamma_d(j-1)\right) - \gamma_d(0).
\]
A moment based estimator of the variance is hence
\begin{equation}
\label{varestimate}
\hat{\sigma}^2=\left(\sum_{j=1}^p \hat{\phi}_j \hat{\gamma}_d(j-1)\right)-\hat{\gamma}_d(0).
\end{equation}
The next section shows that $\hat{\sigma}^2$ is a consistent estimator of $\sigma^2$. 

\section{Asymptotic Normality}

This section shows that if $m=m(N)$ grows slowly enough in $N$, the estimators in the last section will be consistent and asymptotically normal. In particular to obtain asymptotic normality, we will assume that as $n \to \infty$:
\begin{itemize}
    \item[A.1] $\max_{0 \leq k\leq m(n)} |\mu_{k+1}-\mu_k|\leq B$ (the maximum change size is bounded),
    \item[A.2] $m(n)=O(\sqrt{n})$.
\end{itemize}
We begin with asymptotic normality of the autocorrelations for first differences in the general ARMA($p,q$) case, which may be of distinct interest.  The asymptotic normality of the AR($p$) estimators is a corollary to Theorem \ref{Theorem1}.

\begin{thm}
\label{Theorem1}
If $\{ X_t \}_{t=1}^N$ follows (\ref{shiftmodel}) with $\{ \epsilon_t \}$ satisfying (\ref{ARp}), then for each fixed positive integer $k$, as $N \to \infty$,  
\[
\sqrt{N}
\left( 
\begin{array} {c}
\hat{\rho}_d(1) - \rho_d(1) \\
\vdots \\
\hat{\rho}_d(k) - {\rho}_d(k) \\
\end{array}
\right)
\xlongrightarrow {\mathbf{D}} 
\mathcal{N}_k{(\mathbf{0}}, { \mathbf{BWB}^T}).
\]
Here, the elements in $\mathbf{W}$ are given by Bartlett's formula, (see Chapter 8 of \citet{BrockwellDavis1991}) and $\mathbf{B}$ has form
\[    
\mathbf{B}=
\frac{1}{2(1-\rho(1))}
\begin{bmatrix}
2      & -1    & 0     & 0      &\cdots & 0      \\
-1     & 2     &-1     & 0      &\cdots & 0      \\
0      & -1    & 2     &-1      &\cdots & 0      \\
\vdots &\vdots &\vdots & \vdots &\ddots & \vdots \\
0      &0      &0      &0       &\cdots &-1      \\
    \end{bmatrix} .
\]
\end{thm}

\noindent \textbf{Proof} We first show that the changepoints have negligible impact on estimated autocorrelations in the limit. To do this, write $d_t=X_t-X_{t-1}= (\epsilon_t - \epsilon_{t-1}) +\delta_t$, with $\delta_t=(\mu_k-\mu_{k-1}) I_{t=\tau_{k+1}}$.
If we let 
\[
\tilde{\gamma}_d(h)=
\frac{\sum_{t=2}^{N-h} (\epsilon_t-\epsilon_{t-1})(\epsilon_{t+h}-\epsilon_{t+h-1})}{N},
\]
then 
\begin{eqnarray*}
\sqrt{N} \left|\hat{\gamma}_d(h)-\tilde{\gamma}_d(h) \right|
&\leq &\frac{m}{\sqrt{N}}
\left[m^{-1}  \sum_{t=\tau_j}\delta_t(\epsilon_{t+h}-\epsilon_{t+h-1}
+\epsilon_{t-h}-\epsilon_{t-h-1}+\delta_t)\right].
\end{eqnarray*}
The term on the right hand side converges to zero if $N^{-1/2}m \to 0$ as $N \to \infty$ and the sum is bounded in probability.  The assumptions at the start of this section are sufficient to ensure this convergence.  

Now we find the asymptotic distribution of the autocorrelation of the differences of $\{\epsilon_t\}$. Note that \[
\rho_d(h)=\frac{-\rho(h-1)+2\rho(h)-\rho(h+1)}{2(1-\rho(1))}
\quad \text{and} \quad
\hat{\rho}_d(h)=\frac{-\hat{\rho}(h-1)+2\hat{\rho}(h)-\hat{\rho}(h+1)}{2(1-\hat{\rho}(1))}.
\]
Theorem \ref{Theorem1} now follows by the mapping theorem and well known results for asymptotic normality for sample autocovariances for ARMA processes (see Chapter 8 of \cite{BrockwellDavis1991}).  

\begin{col}
\label{armaFitNull}
Suppose that $\{ X_t \}$ follows (\ref{shiftmodel}) with  $\{ \epsilon_t \}$ satisfying (\ref{ARp}) with $q=0$. For the estimator given in (\ref{YWestimator}), as $N \to \infty$,  
\[
\sqrt{N}
\left(
\begin{array}{c}
\hat{\phi}_1 - \phi_1 \\ 
\vdots                \\
\hat{\phi}_p - \phi_p \\
\end{array}
\right)
\xlongrightarrow{\mathbf{D}} 
\mathcal{N}_p({\bf 0}, \boldsymbol{\Sigma}),
\]
Here, $\boldsymbol{\Sigma}={\bf M B W M B}^T$.
\end{col}

\noindent \textbf{Proof of Corollary \ref{armaFitNull}}. Since $\{ d_t \}$ is stationary and ergodic, the elements of $\hat{\boldsymbol{M}}$ converge to those of $\boldsymbol{M}$ in the almost sure sense.  The conclusion of Corollary \ref{armaFitNull} now follows by application of the continuous mapping theorem.

\section{A Simulation Study}

A simulation study with AR$(p)$ errors is now conducted. In each simulation, $N=1000$ and $m$ is randomly generated by the discrete uniform distribution $\mbox{Uniform} \{ 0, 1, \cdots, 10 \}$, which roughly corresponds to the level of changepoint frequency in our ensuing data example.  Any changepoint times are generated randomly within $\{ 2, 3, \ldots, N \}$ with equal probability --- we do not impose any minimal spacing between successive changepoint times. The segment means $\mu_i$ are randomly generated from a $\mbox{Uniform}(-1.5,1.5)$ distribution. Ten thousand independent runs are conducted in all cases.

We start with AR$(1)$ errors, simulating $\phi$ randomly from $\mbox{Uniform}(-0.95, 0.95)$ and $\{ Z_t \}$ is Gaussian white noise with a unit variance. The estimator in (\ref{ar1robust}) is denoted by AR1seg and the Yule-Walker difference estimator in (\ref{YWestimator}) is denoted by Diff.  Raw rolling window estimators, using different window lengths are denoted by their lengths: N, N/2, N/5, N/10, N/20 and N/50.   

\begin{figure}
\centering
\includegraphics[scale=0.4]{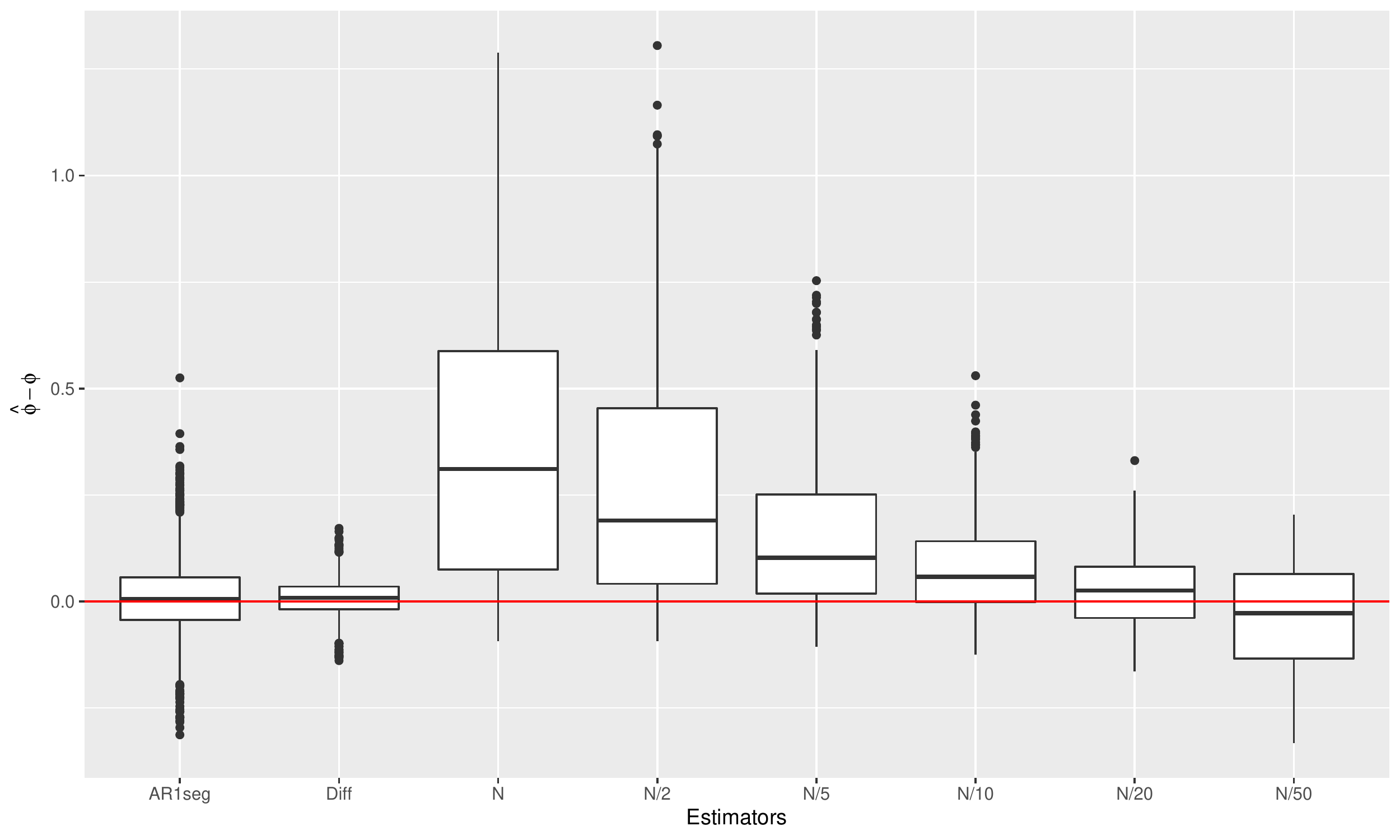}
\caption{Box Plots of Estimates of the AR(1) Parameter $\phi$}
\label{ar1_sim}
\end{figure}

Our simulation results are summarized in Figure \ref{ar1_sim}. The obvious winner is the Yule-Walker estimator based on first order differencing which is unbiased and has a smaller variance. The AR1seg estimator is unbiased as claimed; however, it has larger variability than the Yule-Walker difference estimators. The performance of the raw rolling-window estimators depends on the choice of the window length, but appears to be inferior to the difference based estimator, even with the optimal window size (which is likely somewhere between $N/20$ and $N/50$ in this simulation).  It is hard to decide the optimal window length in practice and smaller window lengths considerably increase computation time.

We now move to AR(2) errors. In each AR(2) run, the AR coefficients were randomly generated (uniformly) from the triangular region that guarantees causality: $\phi_1+\phi_2 < 1$, $\phi_2 - \phi_1 < 1$, and $|\phi_2| < 1$.  In this setting, the changepoint total is fixed at $m=9$ and all segments are set to have equal length.   All mean shifts alternate in sign with an absolute magnitude of 2, the first shift moving the series upwards.   The series length varies from $N=1000$ to $N=20000$. Since $p > 1$, the AR1seg estimator is not applicable. The rolling-window estimator is dropped from consideration due to its poor AR(1) performance and computational time. The simulation results show that estimator bias and variance decreases towards zero as the length of the series increases, reinforcing the consistency result proven in the last section.  

\begin{figure}
\centering
\includegraphics[scale=0.65]{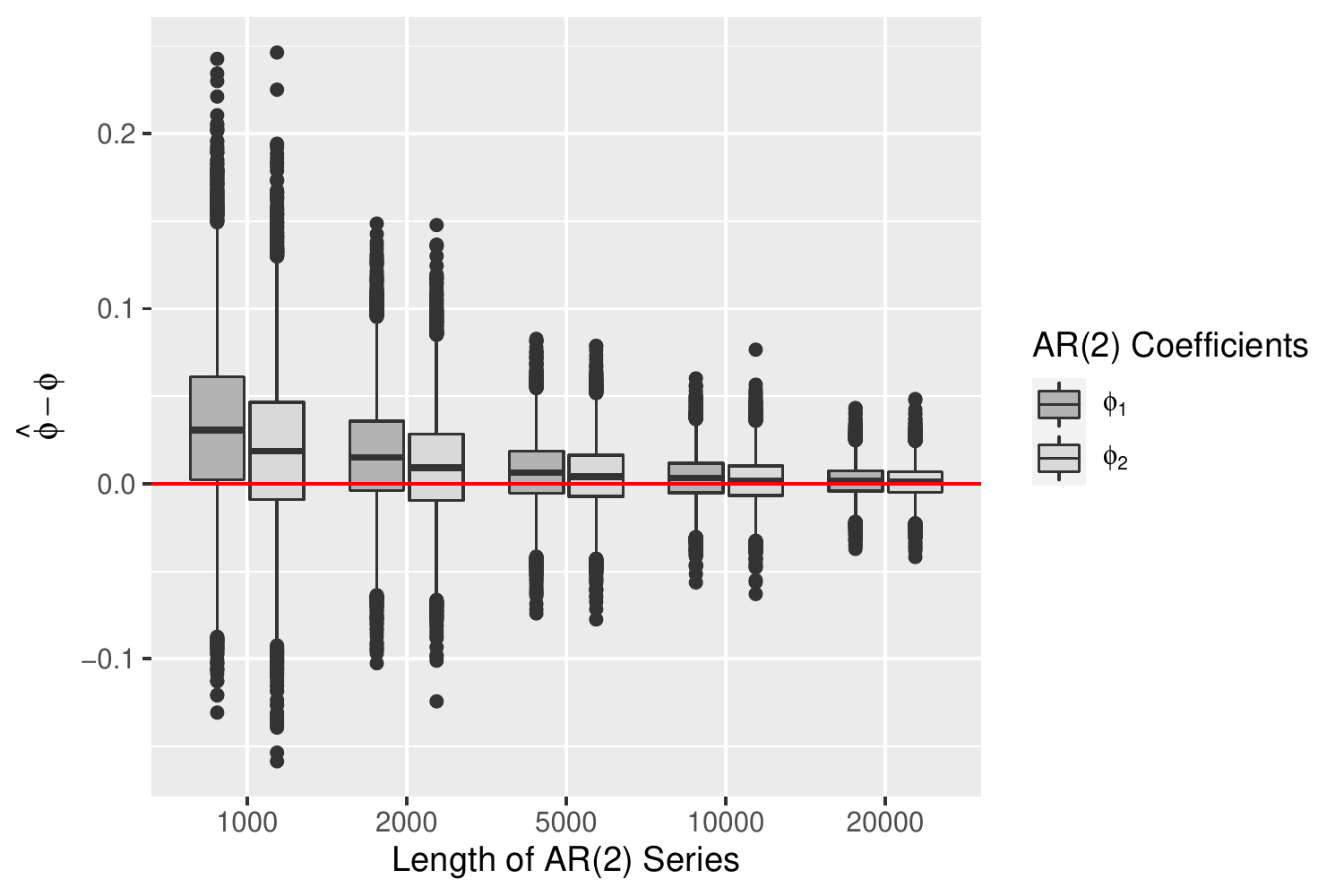}
\caption{Box Plots of AR(2) Coefficient Estimates}
\label{ar2_diff_N}
\end{figure}

\begin{figure}
\centering
\includegraphics[scale=0.65]{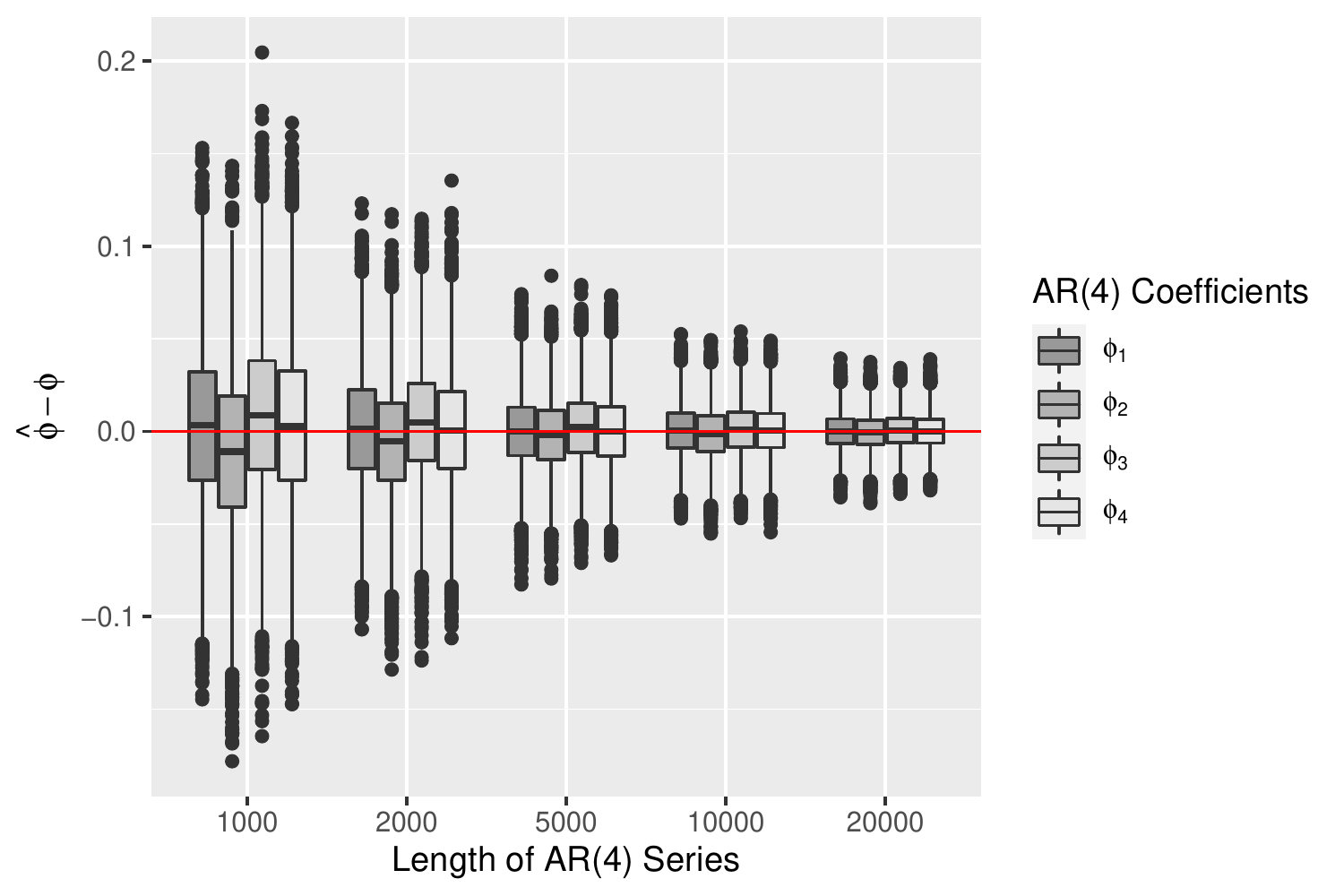}
\caption{Box Plots of AR(4) Coefficient Estimates}
\label{ar4_diff_N}
\end{figure}

Moving to AR$(4)$ simulations, to meet causality requirements, the AR(4) characteristic polynomial is factored into its four roots $1/r_1, 1/r_2, 1/r_3$, and $1/r_4$:
\begin{align*}
    \phi(z)=(1-r_1z)(1-r_2z)(1-r_3z)(1-r_4z).
\end{align*}
Causality implies that all $r_i$ should lie outside the complex unit circle. To meet this, $r_1$ and $r_2$ will be randomly generated from the Uniform$(-0.9, 0.9)$, and $r_3$ is a randomly generated complex number with modulus $\lVert r_3 \rVert<0.9$. $r_4$ is taken as the complex conjugate of $r_3$. This mixes real and complex roots in the AR polynomial.  Other simulation settings are identical to those in the AR$(2)$ case. Figure \ref{ar4_diff_N} shows our results, which exhibit the same pattern in the AR$(2)$ case, with decreasing bias and variance as $N$ increases.

\begin{figure}
\centering
\includegraphics[scale=0.45]{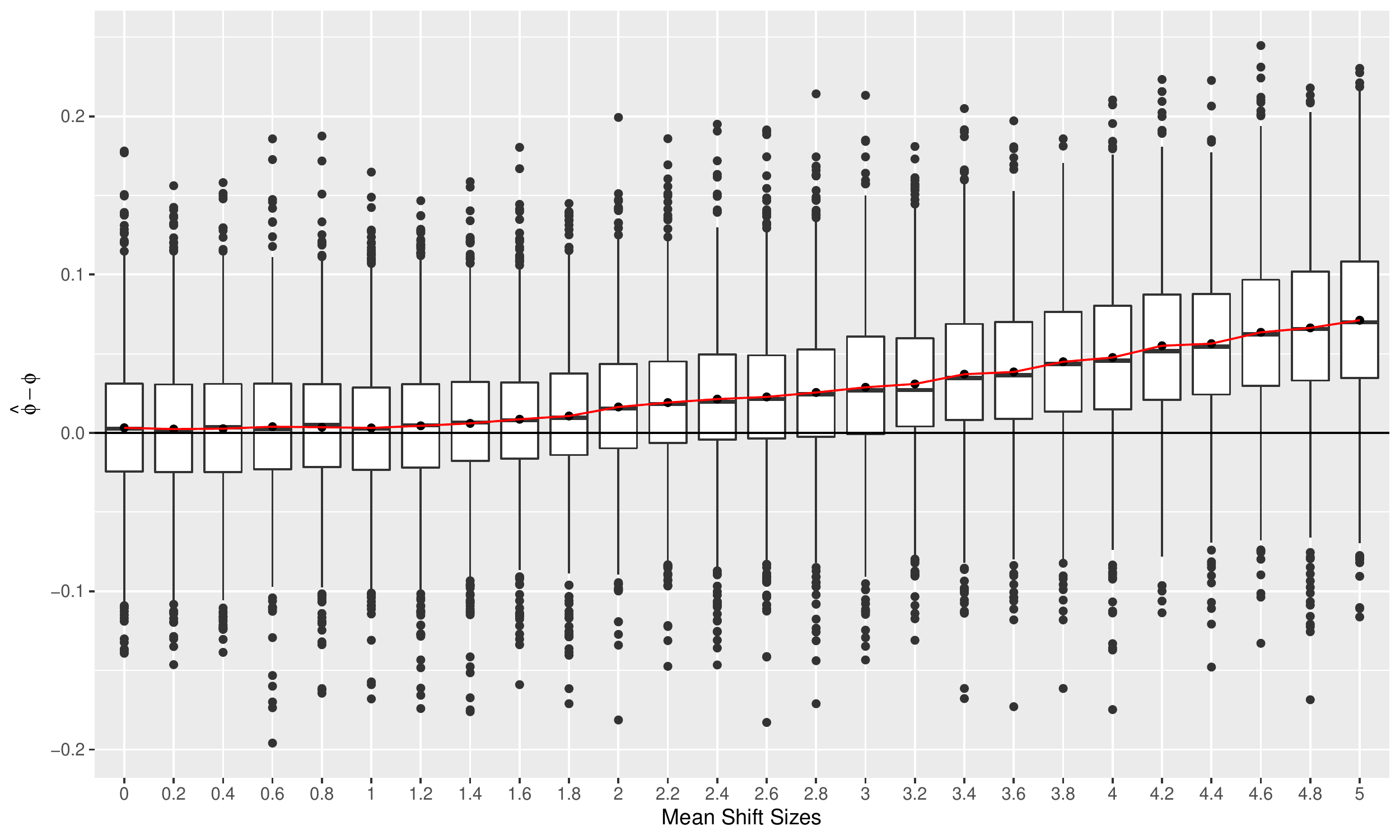}
\caption{An AR(1) Mean Shift Size Sensitivity Plot}
\label{sensitivity}
\end{figure}

For our final simulation, we return to the AR(1) setting and conduct a sensitivity analysis to the mean shift size. When $m/N \rightarrow 0$ as $N \rightarrow \infty$, the accuracy of this estimator is not influenced by changepoint locations but rather the magnitude of the mean shifts --- the signal-to-noise ratio, defined as the absolute mean shift magnitude over the standard deviation.  For simplicity, $\sigma^2$ is set to unity. The number of changepoints is fixed at $m=9$ and their locations randomly generated over $\{ 2, \ldots N \}$ with $N=1000$. In each run, the true $\phi$ is simulated from the  $\text{Uniform}(-0.95, 0.95)$ distribution.  The nine mean shifts are alternating, with fixed shift sizes between zero to $5$ considered across runs. The difference between the estimated $\hat{\phi}$ and the true $\phi$ is presented 
in Figure \ref {sensitivity}, aggregated by mean shift size.

The horizontal line in Figure \ref{sensitivity} demarcates zero bias in $\hat{\phi}$; the red curve depicts the average differences between $\hat{\phi}$ and $\phi$. Obviously, the larger the mean shift magnitude, the more our estimator degrades.  However, larger mean shift sizes can easily be identified as outliers in the differenced series \citep{outlier1993, outlier2003}.  As such, the essential challenge lies with estimating the AR($p$) parameters in the presence of the smaller mean shifts.

\section{Applications}

\subsection{Changepoints in AR(\texorpdfstring{$p$}) Series}

Some techniques can mistakenly flag changepoints when underlying positive dependence is ignored. For example, \citet{lund2020wbs} argues that the London house price series shifts identified may be more attributable to the positive correlations in the series than to actual mean shifts. CUSUM based techniques are severely affected by positive correlation \cite{shi2021comparison}.  To remedy this, authors recommend detecting changepoints from estimated versions of the one-step-ahead prediction residuals of the series \citep{bai1993resid, robbins2011}. However, this requires estimation of the autocovariance structure of the series in the presence of the unknown changepoints.  A major application of our methods serves to decorrelate series without requiring any knowledge about the changepoint configuration.  IID-based changepoint techniques, applied to the estimated one-step-ahead prediction resdiuals, can then locate mean shifts in the series. The Yule-Walker difference estimator proposed here is extensively used in \cite{shi2021comparison} to do just this.

Table \ref{app1} illustrates the improved performance of two popular multiple changepoint methods, Wild Binary Segmentation (WBS) \citep{fryzlewicz2014wbs} and Pruned Exact Linear Time  (PELT) \citep{killick2012pelt}. In each run, an AR$(1)$ series of length $N=500$ is simulated with $\phi$ fixed in $\{ 0.25, 0.50, 0.75 \}$ and white noise variance $\sigma^2=1$.  The series has either zero or three equally spaced changepoints; the mean shift sizes $\Delta$ are chosen to satisfy the constant signal-to-noise ratio requirement 
\[
\text{SNR} = \frac{|\Delta|} {\sqrt{\frac{\sigma^2}{1-\phi^2}}}=2. 
\]
All simulations employ $1000$ independent runs. In Table  \ref{app1}, $\overline{\hat{m}}$ and $SE_{\hat{m}}$ denote the average and standard error of the estimated number of changepoints when WBS and PELT are directly applied to the series.  The quantities $\overline{\hat{m}^D}$ and $SE_{\hat{m}^D}$ denote the average and standard error of the estimated number of changepoints from the one-step-ahead prediction residuals. 

\begin{table}
\centering
\begin{tabular}{c}
  \begin{tabular}{c|cc|cc}
    \multicolumn{5}{c}{$\phi=0.25$}\\
    \hline
    \multirow{2}*{\# of Changepoints} &
      \multicolumn{2}{c}{WBS} &
      \multicolumn{2}{c}{PELT} \\
      & {$\overline{\hat{m}}$/($SE_{\hat{m}}$)} & {$\overline{\hat{m}^D}$/($SE_{\hat{m}^D}$)} & {$\overline{\hat{m}}$/($SE_{\hat{m}}$)} & {$\overline{\hat{m}^D}$/($SE_{\hat{m}^D}$)} \\
    \hline
    Zero & 3.85/(2.43) & 0.17/(0.54)  & 0.02/(0.17) & 0.00/(0.02) \\
    Three & 5.46/(2.01) &  3.03/(0.18) & 3.07/(0.34) & 3.00/(0.03) \\
    \hline
  \end{tabular} \\
  \begin{tabular}{c|cc|cc}
    \multicolumn{5}{c}{$\phi=0.5$}\\
    \hline
    \multirow{2}*{\# of Changepoints} &
      \multicolumn{2}{c}{WBS} &
      \multicolumn{2}{c}{PELT} \\
      & {$\overline{\hat{m}}$/($SE_{\hat{m}}$)} & {$\overline{\hat{m}^D}$/($SE_{\hat{m}^D}$)} & {$\overline{\hat{m}}$/($SE_{\hat{m}}$)} & {$\overline{\hat{m}^D}$/($SE_{\hat{m}^D}$)} \\
    \hline
    Zero & 16.03/(3.92) & 0.24/(0.70)  & 1.36/(1.62) & 0.00/(0.04) \\
    Three & 16.90/(3.81) & 3.09/(0.40) & 4.28/(1.43) & 2.95/(0.37) \\
    \hline
  \end{tabular}   \\
  \begin{tabular}{c|cc|cc}
    \multicolumn{5}{c}{$\phi=0.75$}\\
    \hline
    \multirow{-0.65}*{\# of Changepoints} &
      \multicolumn{2}{c}{WBS} &
      \multicolumn{2}{c}{PELT} \\
      & {$\overline{\hat{m}}$/($SE_{\hat{m}}$)} & {$\overline{\hat{m}^D}$/($SE_{\hat{m}^D}$)} & {$\overline{\hat{m}}$/($SE_{\hat{m}}$)} & {$\overline{\hat{m}^D}$/($SE_{\hat{m}^D}$)} \\
    \hline
    Zero & 30.77/(3.98) & 0.40/(0.96)  & 12.65/(3.31) & 0.01/(0.13) \\
    Three & 31.45/(3.86) &  2.48/(1.38) & 14.33/(3.25) & 1.59/(1.44) \\
    \hline
  \end{tabular} 
\end{tabular}  
\caption{AR$(1)$ series with/without changepoints. The thresholding constant used in WBS is $C=1.3$.}
\label{app1}
\end{table}

It is apparent that IID based WBS and PELT methods overestimate the number of changepoints in a dependent series when correlation is ignored; PELT seems more resistant to dependence than WBS.  In contrast, with the help of the proposed Yule-Walker difference estimator and decorrelation techniques, both WBS and PELT become much more accurate.

\subsection{New Bedford Precipitation}
\label{NBprecip}

\begin{figure}
\centering
\includegraphics[scale=0.6]{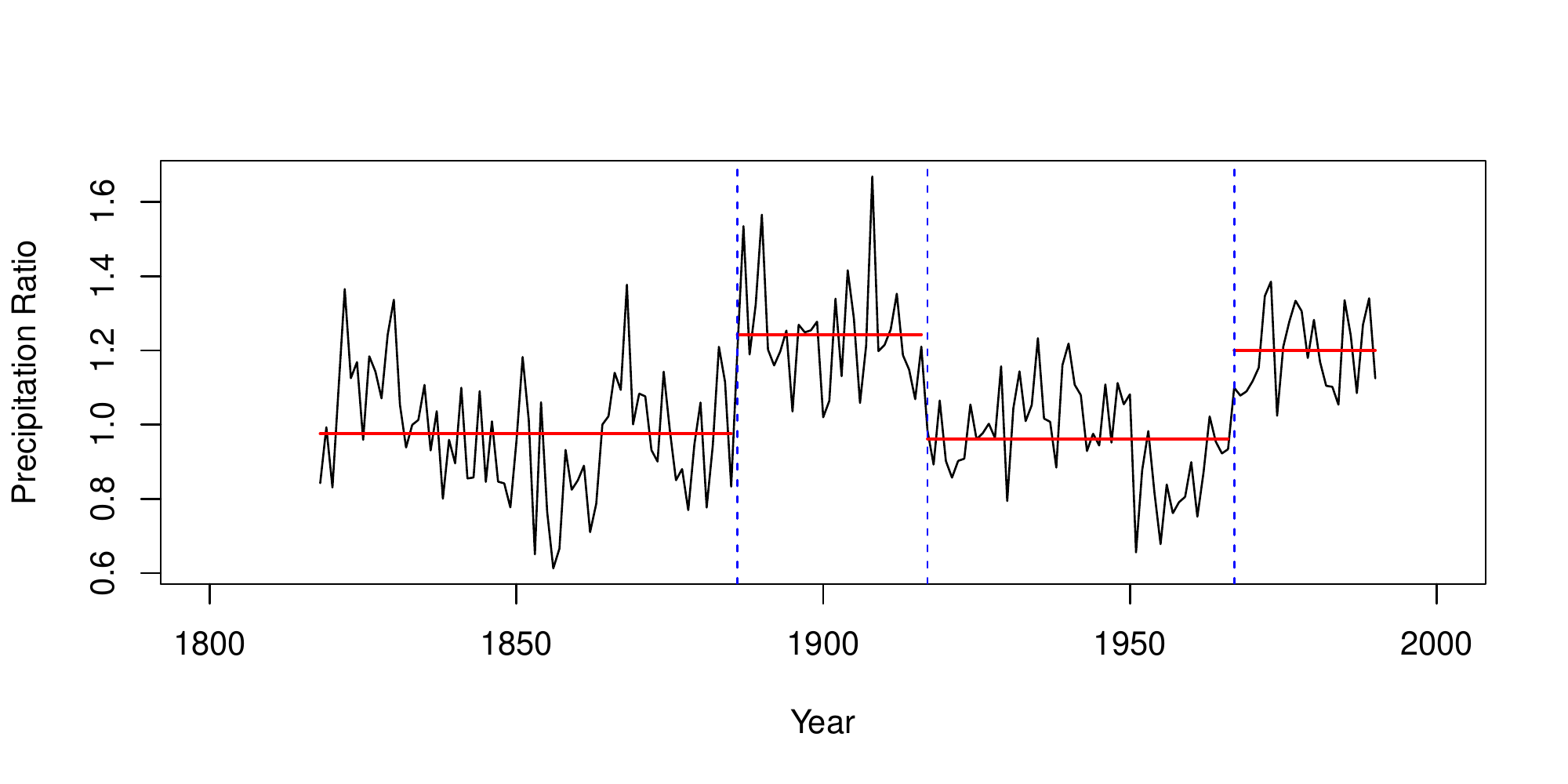}
\caption{New Bedford to Boston Annual Precipitation Ratios with Three Changepoints Demarcated.}
\label{precipitationplot}
\end{figure}

\cite{lilund2012MDL} studied annual precipitations from New Bedford and Boston, Massachusetts. The data are available from \href{https://w2.weather.gov/climate/xmacis.php?wfo=box}{the National Oceanic and Atmospheric Administration (NOAA)}. The ratio of these series (New Bedford to Boston) is displayed in Figure \ref{precipitationplot} along with a fitted mean of a model that allows for both multiple changepoints and AR(1) errors. Three changepoints, occurring at the years 1886, 1917, and 1967, are indicated. 

What is interesting here is the AR(1) parameter estimate: it fluctuates wildly over distinct methods.  Specifically, our difference estimator and AR(1)seg methods produce antipodal estimates, as can be seen in Table \ref{precipitationtab}. Our estimate agrees closely with an estimate computed by assuming the three changepoint times are known, but the level of autocorrelation is significantly less than that estimated in a Yule-Walker scheme that ignores all three changepoint times.   Overall, the results show that one needs to be careful with changepoint problems with correlated data --- mean shifts and correlation can inject similar features into time series.

\begin{table}
\begin{center}
    \begin{tabular}{l|c}
         Methods & Estimate ($\hat{\phi}$)  \\
             \hline
         Yule Walker Estimator (ignoring changepoints) & $0.547$ \\   
         Robust Estimator (AR1seg) & $-0.268$\\
         Difference Estimator & $0.255$\\
         Yule Walker Estimator (given changepoint times) & $0.273$ 
    \end{tabular}
    \caption{AR(1) Coefficient Estimation}
    \label{precipitationtab}
\end{center}
\end{table}

\section{Conclusion}
Differencing can help estimate the autocovariance structure of an AR$(p)$ series corrupted by mean shifts. Our Yule-Walker estimator for autoregressive models is easy to implement, computationally fast, consistent, and asymptotically normal.  The proposed estimator is adversely impacted by large mean shifts; however these large shifts appear as outliers in the differenced series and can easily be removed.  With changepoints present, the difference estimator significantly improves changepoint techniques developed for models with IID errors.  In addition, the techniques are applicable if the series has a linear trend (constant across all regimes) with intercept shifts.

\bibliographystyle{plainnat}

\begin{thebibliography}{}

\bibitem[Bai, 1993]{bai1993resid}
Bai, J. (1993).
\newblock On the partial sums of residuals in autoregressive and moving average
  models.
\newblock {\em Journal of Time Series Analysis}, 14(3):247--260.

\bibitem[Beaulieu and Killick, 2018]{BeaulieuKillick2018}
Beaulieu, C. and Killick, R. (2018).
\newblock Distinguishing trends and shifts from memory in climate data.
\newblock {\em Journal of Climate}, 31:9519--9543.

\bibitem[Brockwell and Davis, 1991]{BrockwellDavis1991}
Brockwell, P.~J. and Davis, R.~A. (1991).
\newblock {\em Time Series: Theory and Methods}.
\newblock Springer Series in Statistics. Springer-Verlag, New York, second
  edition.

\bibitem[Chakar et~al., 2017]{chakar2017}
Chakar, S., Lebarbier, E., L{\'e}vy-Leduc, C., Robin, S., et~al. (2017).
\newblock A robust approach for estimating change-points in the mean of an
  $\text{AR}(1)$ process.
\newblock {\em Bernoulli}, 23(2):1408--1447.

\bibitem[Chen and Liu, 1993]{outlier1993}
Chen, C. and Liu, L.-M. (1993).
\newblock Joint estimation of model parameters and outlier effects in time
  series.
\newblock {\em Journal of the American Statistical Association},
  88(421):284--297.

\bibitem[Fryzlewicz, 2014]{fryzlewicz2014wbs}
Fryzlewicz, P. (2014).
\newblock Wild binary segmentation for multiple change-point detection.
\newblock {\em The Annals of Statistics}, 42(6):2243--2281.

\bibitem[Hall et~al., 2000]{hall2000disease}
Hall, C.~B., Lipton, R.~B., Sliwinski, M., and Stewart, W.~F. (2000).
\newblock A change point model for estimating the onset of cognitive decline in
  preclinical alzheimer's disease.
\newblock {\em Statistics in Medicine}, 19(11-12):1555--1566.

\bibitem[Hewaarachchi et~al., 2017]{lund2017temp}
Hewaarachchi, A.~P., Li, Y., Lund, R., and Rennie, J. (2017).
\newblock Homogenization of daily temperature data.
\newblock {\em Journal of Climate}, 30(3):985--999.

\bibitem[Killick et~al., 2012]{killick2012pelt}
Killick, R., Fearnhead, P., and Eckley, I.~A. (2012).
\newblock Optimal detection of changepoints with a linear computational cost.
\newblock {\em Journal of the American Statistical Association},
  107(500):1590--1598.

\bibitem[Li and Lund, 2012]{lilund2012MDL}
Li, S. and Lund, R. (2012).
\newblock Multiple changepoint detection via genetic algorithms.
\newblock {\em Journal of Climate}, 25(2):674--686.

\bibitem[Lund and Shi, 2020]{lund2020wbs}
Lund, R. and Shi, X. (2020).
\newblock Commentary on: Detecting possibly frequent change-points: Wild binary
  segmentation 2 and steepest-drop model selection.
\newblock {\em Journal of the Korean Statistical Society}, 49:1090--1095.

\bibitem[Maidstone et~al., 2017]{maidstone2017fpop}
Maidstone, R., Hocking, T., Rigaill, G., and Fearnhead, P. (2017).
\newblock On optimal multiple changepoint algorithms for large data.
\newblock {\em Statistics and Computing}, 27(2):519--533.

\bibitem[Marriott and Pope, 1954]{marriott1954}
Marriott, F. and Pope, J. (1954).
\newblock Bias in the estimation of autocorrelations.
\newblock {\em Biometrika}, 41(3/4):390--402.

\bibitem[McQuarrie and Tsai, 2003]{outlier2003}
McQuarrie, A.~D. and Tsai, C.-L. (2003).
\newblock Outlier detections in autoregressive models.
\newblock {\em Journal of Computational and Graphical Statistics},
  12(2):450--471.

\bibitem[Norwood and Killick, 2018]{NorwoodKillick2018}
Norwood, B. and Killick, R. (2018).
\newblock Long memory and changepoint models: A spectral classification
  procedure.
\newblock {\em Statistics \& Computing}, 28(2):291--302.

\bibitem[Orcutt and Winokur~Jr., 1969]{orcutt1969}
Orcutt, G.~H. and Winokur~Jr., H.~S. (1969).
\newblock First order autoregression: inference, estimation, and prediction.
\newblock {\em Econometrica: Journal of the Econometric Society}, pages 1--14.

\bibitem[Page, 1954]{Page1954}
Page, E.~S. (1954).
\newblock Continuous inspection schemes.
\newblock {\em Biometrika}, 41(1-2):100--115.

\bibitem[Robbins et~al., 2011]{robbins2011}
Robbins, M., Gallagher, C., Lund, R., and Aue, A. (2011).
\newblock Mean shift testing in correlated data.
\newblock {\em Journal of Time Series Analysis}, 32(5):498--511.

\bibitem[Shi et~al., 2021]{shi2021comparison}
Shi, X., Gallagher, C., Lund, R., and Killick, R. (2021).
\newblock A comparison of single and multiple changepoint techniques for time
  series data.
\newblock {\em arXiv preprint arXiv:2101.01960}.

\end{thebibliography}

\end{document}